\documentclass[12pt]{amsart}

\usepackage{latexsym}
\usepackage{amsmath}
\usepackage{amssymb}
\usepackage{amscd}
\usepackage{multicol}

\newtheorem{theorem}{\bf Theorem}[section]

\newtheorem{proposition}[theorem]{\bf Proposition}

\newtheorem{definition-theorem}[theorem]{\bf Theorem-Definition}

\def\bR{\mathbb{R}}
\def\bC{\mathbb{C}}
\def\bZ{\mathbb{Z}}

\def\quott({/\! /}
\def\la{(\! (}
\def\ra{)\! )}

\def\p{\partial}

\def\fY{{\mathcal Y}}

\def\V{{\mathcal V}}
\def\X{{ X}}
\def\R{{\mathcal R}}
\def\F{{\mathcal F\ell}^{(n)}}
\def\S{\frak{S}}
\def\Y{{ Y}}

\def\D{{\mathcal D}}

\title[A quantum deformation of the cohomology of flag manifolds]
{A quantum type deformation of the cohomology ring of flag manifolds}
\date{\today}

\date{\today}
\begin{document}
\begin{abstract} Let $q_1,\ldots,q_n$ be some variables and set $K:=\bZ[q_1,\ldots,q_n]/( \prod_{i=1}^n q_i)$. We show that there exists  a $K$-bilinear product $\star$ on $H^*(F_n;\bZ)\otimes K$ which  is uniquely determined by some quantum cohomology like properties (most importantly, a degree two relation involving the generators and  an analogue of the flatness of the Dubrovin connection). Then we prove that 
$\star$ satisfies the Frobenius property with respect to the Poincar\'e pairing of 
$H^*(F_n;\bZ)$; this leads immediately to the orthogonality of the corresponding Schubert type polynomials.  We also 
note that 
if we pick $k\in\{1,\ldots,n\}$ and we formally replace  $q_k$ by 0, the ring $(H^*(F_n;\bZ)\otimes K,\star)$ becomes  isomorphic  to the usual small quantum cohomology ring of $F_n$, by an isomorphism which is described precisely. 
\end{abstract}
\author[A.-L. Mare]{Augustin-Liviu  Mare}
\address{
Department of Mathematics and Statistics\\ University of
Regina \\ College West 307.14 \\ Regina SK, Canada S4S
0A2}
\email{mareal@math.uregina.ca}

\maketitle
\section{Introduction}

We consider the complex flag manifold
$$F_n=\{V_1\subset V_2 \subset  \ldots \subset V_{n-1} \subset \bC^n \ | \  V_k \
{\text  {\rm is  a  }}  k{\text {\rm  -dimensional    linear  subspace   of  }} \bC^n \}.$$
For every $k\in\{1,\ldots,n-1\}$ we consider the tautological vector bundle $\V_k$ over $F_n$ and the cohomology\footnote{The coefficient ring for
cohomology will be always $\bZ$, unless otherwise specified.} class
$y_k=-c_1(\det \V_k)\in H^2(F_n).$
It is known that the cohomology classes $y_1,\ldots,y_{n-1}$ generate the ring $H^*(F_n)$. 
To describe the ideal of relations, it is convenient to consider the classes
 $x_k:=y_{k}-y_{k-1},$ $1\le k \le n$ (where we assign
$y_0=y_n:=0$). If $Y_1,\ldots,Y_{n-1}$ are some variables, we set
$$X_k:=Y_{k}-Y_{k-1},$$
where $Y_0=Y_n:=0$. 
We denote by $S(X_1,\ldots,X_n)$ the ideal of $\bZ[X_1,\ldots,X_n]$
generated by the symmetric polynomials with zero constant term.
By a theorem of  Borel, there is a ring isomorphism
\begin{equation}\label{borel}H^*(F_n) \simeq \bZ[Y_1,\ldots,Y_{n-1}]/ 
S(X_1,\ldots,X_n).\end{equation}
Here $y_i$ is mapped to the coset of $Y_i$, for all $i\in\{1,\ldots,n-1\}$.

The small quantum cohomology ring $qH^*(F_n)$ is a  deformation with
$n-1$ parameters of  $H^*(F_n)$. A theorem of Givental and Kim
\cite{Gi-Ki} describes $qH^*(F_n)$ as a  quotient of a certain polynomial ring by the ideal generated by the integrals of motions of the open Toda lattice. 
In this paper we consider the ``periodic" version of the latter (quotient) ring:
that is, the ideal of relations is generated by the integrals of motion of the
{\it periodic} Toda lattice. We address the following question:
does this new ring arise in the spirit of Givental and Kim from a deformation
of $H^*(F_n)$? We construct such a deformation, give a list of properties  
that characterizes it uniquely, and then study it briefly.

Let us be more precise. We  consider $n$ variables $q_1,\ldots,q_n$ and set 
$$K:=\bZ[q_1,\ldots,q_n]/( \prod_{i=1}^nq_i).$$
The main result of this paper is as follows.
\begin{theorem}\label{firstmain} There exists a product $\star$ on $H^*(F_n)\otimes K$ which is uniquely determined by the following properties:
\begin{itemize}
\item[(i)] $\star$ is $K$-bilinear
\item[(ii)] $\star$ preserves the grading induced by the
usual grading of $H^*(F_n)$ combined with $\deg q_j=4$, $1\leq j
\leq n$
\item[(iii)] $\star$ is a deformation of the usual product, in the sense
that if we formally replace all $q_j$ by $0$, we obtain the usual
(cup-)product on $H^*(F_n)$
\item[(iv)] $\star$ is commutative
\item[(v)] $\star$ is associative
\item[(vi)]  we have the  relation 
$$\sum_{1\le i<j \le n}x_i\star x_{j} +\sum_{j=1}^nq_j=0$$
\item[(vii)] the coefficients $(y_i \star a)_d\in H^*(F_n)$ of 
$q^d:=q_1^{d_1}\ldots q_n^{d_n}$ in
$y_i \star a$ satisfy $$(d_i-d_{n})(y_j
\star a)_d=(d_j-d_{n})(y_i \star a)_d,$$ for
 any $d=(d_1,\ldots,d_{n})\geq 0$, $a\in H^*(F_n)$, and $1\leq i,j\leq n-1$.
\end{itemize}
\end{theorem}

The product $\star$ is constructed in Section \ref{three}. 
The main input consists of the conservation laws of the periodic quantum Toda
lattice (cf. \cite{Go-Wa}, \cite{Ko}).
From the construction  we can see that we have a ring isomorphism 
\begin{equation}\label{ssta}(H^*(F_n)\otimes K,\star)\simeq \bZ[Y_1,\ldots,Y_{n-1},q_1,\ldots,q_n]/({\mathcal R}_1,
\ldots,{\mathcal R}_n),\end{equation}
where $ {\mathcal R}_1,
\ldots,{\mathcal R}_n$ are essentially the integrals of motion of the periodic Toda lattice. They
are described  explicitly by equation
(\ref{dete}) below. We note at this point that ${\mathcal R}_n=\pm \prod_{i=1}^n q_i$.
This explains why the ``quantum parameter" ring of our deformation of $H^*(F_n)$ is no longer a polynomial ring, like for $qH^*(F_n)$, but rather  the somewhat awkward  ring  $\bZ[q_1, \ldots, q_n]/ (\prod_{i=1}^n q_i)$. 

We will also prove some properties of the product $\star$.
The first one involves the Poincar\'e pairing $(  \ , \ )$ on $H^*(F_n)$; we actually extend it to a $K$-bilinear form on $H^*(F_n)\otimes K$.
     
\begin{theorem}\label{onemain}  The product $\star$
satisfies the following Frobenius type property: 
$$(a\star b, c)=(a,b\star c)$$
for any $a,b,c\in H^*(F_n)$.
\end{theorem}

We also consider representatives of Schubert cohomology classes 
in $H^*(F_n)$ via the isomorphism given by (\ref{ssta}), which are polynomials
in $Y_1,\ldots,Y_{n-1}, q_1, \ldots, q_n$.
Theorem \ref{onemain} will allow us to prove that  these polynomials satisfy a certain orthogonality relation, similar to the one satisfied by the quantum Schubert polynomials of [Fo-Ge-Po] and [Ki-Ma] (see Proposition \ref{lang} below).

Finally, we consider the ring obtained from $(H^*(F_n)\otimes K,\star)$ by
formally setting a certain $q_k$ to zero. 
More precisely, we denote by $(q_k)$ the set of all elements of $H^*(F_n)\otimes K$
that are multiples of $q_k$ and  take the quotient
\begin{align*}(H^*(F_n)\otimes K)/(q_k) &= 
H^*(F_n)\otimes \bZ[q_1,\ldots,q_{k-1},q_{k+1},\ldots, q_n]\\
{}&=\bZ[x_1,\ldots,x_{n},q_1,\ldots,q_{k-1},q_{k+1},\ldots,q_n]/
S(x_1,\ldots,x_n).
\end{align*}
Here  $S(x_1,\ldots,x_n)$ denotes the ideal  of $\bZ[x_1,\ldots,x_n]$  which is generated by
the symmetric polynomials with zero constant term. The space $(H^*(F_n)\otimes K)/(q_k)$ can be equipped with the product
induced by $\star$. 
The resulting ring is isomorphic to the actual small quantum cohomology ring
$qH^*(F_n)$: the following theorem describes this  ring isomorphism.
We recall (cf. e.g. \cite{Fu-Pa}) that the ring $qH^*(F_n)$  
consists of the space
$$H^*(F_n)\otimes \bZ[Q_1,\ldots,Q_{n-1}]
=\bZ[X_1,\ldots,X_{n},Q_1,\ldots,Q_{n-1}]
/
S( X_1,\ldots, X_{n}),
$$
which is equipped with a certain (quantum) product.

\begin{theorem}\label{secondmain} Fix $k\in\{1,\ldots,n\}$.
The ring $((H^*(F_n)\otimes K)/( q_k),\star)$ is isomorphic to  
$qH^*(F_n)$.
The isomorphism 
is the map from 
$$\bZ[X_1,\ldots,X_{n},Q_1,\ldots,Q_{n-1}]
/
S( X_1,\ldots, X_{n})$$
to
$$
\bZ[x_1,\ldots,x_{n},q_1,\ldots,q_{k-1},q_{k+1},\ldots,q_n]/
S(x_1,\ldots,x_n)$$
which is given by
$$X_i \mapsto x_{k-i+1}, 1\le i \le n, \quad 
Q_j \mapsto q_{k-j}, 1\le j \le n-1,$$
where the indices are evaluated modulo $n$.   
\end{theorem}

\noindent 
{\bf Remarks.} 
1. A characterization of the quantum cohomology ring $qH^*(F_n)$,
similar in spirit to the one given by Theorem \ref{firstmain}, has been obtained 
 by the author in \cite{Ma2}; see also Theorem \ref{mar} below.
 
2. The product $\star$ considered in this paper could also be relevant in the context of
the small quantum cohomology  of the infinite dimensional flag manifold 
$\F$.
  The ring $qH^*(\F)$ is  conjecturally defined and  investigated 
  by Guest and Otofuji in \cite{Gu-Ot} and the author of the present paper in 
  \cite{Ma1}.
  More precisely, in the latter paper
   we have considered  
a product  $\circ$ on $H^*(\F;\bR)\otimes \bR[q_1,\ldots,q_n]$  which satisfies certain natural properties  similar to (i)-(vii)  
in Theorem \ref{firstmain} above; we have also explained why a quantum product
$\circ$
resulting from Gromov-Witten invariants of $\F$ should satisfy these properties.
However, it is still an open question whether such a product  exists. 
As long as nobody gave a rigorous definition (in terms of
stable curves) of $qH^*(\F)$, it seems natural to attempt to
construct it by other means, for instance with the methods of this paper.
In fact, the ring we are constructing and studying here is essentially the quotient of $qH^*(\F)$   
by the product $\prod_{i=1}^nq_i$.

3. 
Let $G$ be a compact simple simply connected 
Lie group and $T\subset G$ a maximal torus. It is likely that some of the results of this paper can be extended to 
the flag manifold $G/T$, where $G$ is of type other than $A$.  
For instance,  if $\ell :=\dim T$, we expect that an extension of $H^*(G/T;\bR)$ with
$\ell +1$ parameters can be constructed with the methods used
in Section \ref{three} below. More precisely, the  ring of  ``quantum parameters" would be in this case
$K=\bR[q_1,\ldots,q_{\ell +1}]/(q_1^{m_1}\cdots q_{\ell}^{m_{\ell}}q_{\ell+1})$.
The numbers $m_1,\ldots,m_\ell$ arise from the expansion
$$\alpha_0^\vee=m_1\alpha_1^\vee +\ldots + m_\ell\alpha_\ell^\vee,$$
where $\{\alpha_1,\ldots,\alpha_\ell\}$ is a simple root system of $G$,
$\alpha_0$ the highest root, and the superscript $\vee$ indicates the 
corresponding coroot. 
However, 
it should be noted that the construction could only work for
$G$ of certain types (e.g. $A$, $B$ or $C$), 
 for exactly the same reason as in \cite{Ma1}. We will probably explore this more
 general situation elsewhere.

 The paper is organized as follows: First we prove Theorem \ref{secondmain} and the uniqueness part of Theorem \ref{firstmain}. After that we prove the existence part of the 
 latter theorem, which uses a presentation of the ring $(H^*(F_n)\otimes K,\star)$ in terms of generators and relations. Finally we prove the Frobenius property and 
the orthogonality of the Schubert type polynomials. Such polynomials are 
determined explicitly in the case $n=3$.

 \noindent {\bf Acknowledgements.} I would like to thank Martin Guest and Takashi Otofuji for discussions about the topics of the paper. 
 I am also extremely grateful to the Referee for suggesting numerous improvements.
 
\section{Uniqueness of the product $\star$}
The goals of this section are:
 show that there exists {\it at most} one product $\star$ with the properties
(i)-(vii) in Theorem \ref{firstmain}; prove Theorem \ref{secondmain}. 
The main instrument is the following result, which is a particular case of
 \cite[Theorem 1.1]{Ma2}.

\begin{theorem}\label{mar}{\rm ([Ma2])} Let $Y_1,\ldots,Y_{n-1},Q_1,\ldots,Q_{n-1}$ be some
variables and  set
$X_j:=Y_{j}-Y_{j-1}$  for $1\le j\le n$, where $Y_0=Y_{n}:=0$. Denote by
$S(X_1,\ldots,X_n)$ the ideal of $\bZ[Y_1,\ldots,Y_{n-1}]$ generated by the  symmetric  polynomials 
with zero constant term in $X_1,\ldots, X_{n}$. There exists a product $\circ$  on
$$ \left(\bZ[Y_1,\ldots,Y_{n-1}]/S(X_1,\ldots,X_n)\right) \otimes \bZ[Q_1,\ldots,Q_{n-1}]$$
which is uniquely determined by the following properties:
\begin{itemize}
\item[(i)] $\circ$ is $\bZ[Q_1,\ldots, Q_{n-1}]$-bilinear
\item[(ii)] $\deg ([f]\circ[g])=\deg[f]+\deg[g]$, for any 
$f,g\in \bZ[Y_1,\ldots,Y_{n-1}]$ which are homogeneous (the brackets $[ \ ]$ indicate the coset modulo $S(X_1,\ldots,X_n)$ and the grading is given by $\deg Y_i:=2$, $\deg Q_i:=4$)
\item[(iii)] $\circ$  is a deformation of the canonical product on $\bZ[Y_1,\ldots,Y_{n-1}]/S(X_1,\ldots,X_n)$
\item[(iv)] $\circ$ is 
commutative
\item[(v)] $\circ$ is associative
\item[(vi)] we have
$$\sum_{1\le i <j \le n}[X_i]\circ [X_j] +\sum_{i=1}^{n-1}Q_i=0$$
\item[(vii)] we have
$$Q_j\frac{\partial}{\partial Q_j}([Y_i]\circ [f])=Q_i\frac{\partial}{\partial Q_i}([Y_j]\circ [f])$$
for any $1\le i,j\le n-1$ and any $f\in \bZ[Y_1,\ldots,Y_{n-1}]$.
\end{itemize}
\end{theorem}  

The product $\circ$ mentioned in this theorem is induced by the small quantum cohomology
ring $qH^*(F_n)$ via 
the Borel isomorphism (\ref{borel}). 

Let now $\star$  be a product which satisfies the assumptions (i)-(vii) in Theorem \ref{firstmain}. The ring 
$(H^*(F_n)\otimes K,\star)$ is generated by $y_1,\ldots,y_{n-1},q_1,\ldots,q_n.$ Consequently, as a $K$-algebra, it is generated by $y_1,\ldots,y_{n-1}$. Let us consider the Schubert basis 
$\{\sigma_w\ | \ w\in S_n\}$ of $H^*(F_n)$.
Here $S_n$ is the symmetric group and $\sigma_w$ are cosets in
$\bZ[x_1,\ldots, x_{n}]/S(x_1,\ldots,x_n)$ of the Schubert polynomials
(cf. e.g. \cite[Ch. 10]{Fu}).
 We also consider the expansion
\begin{equation}\label{yi}y_i\star \sigma_w =\sum_{v\in S_n} \omega_i^{vw}\sigma_v
\end{equation}
where $\omega_i^{vw}\in K$. 
We will show that for any $i$, $v$, and $w$, 
the coefficient $\omega_i^{vw}$ is prescribed.

To this end, we first consider the quotient 
$$(H^*(F_n)\otimes K)/(q_n)
= H^*(F_n)\otimes \bZ[q_1,\ldots, q_{n-1}]=
\bZ[x_1,\ldots,x_{n},q_1,\ldots,q_{n-1}]/
S(x_1,\ldots,x_n)$$
and equip it with the product induced by $\star$. Let us denote this new product by $\star_n$.
It is  commutative, associative, and  satisfies the obvious grading condition.
Moreover, for any $a\in H^*(F_n)$ and any $d=(d_1,\ldots, d_{n-1},0)\ge 0$ we have
$$d_i(y_j\star_n a)_d=d_j(y_i\star_n a)_d,$$
for all $1\le i,j\le n-1$. This implies that 
$$q_i\frac{\partial}{\partial q_i}(y_j\star_n a) =  q_j\frac{\partial}{\partial q_j}(y_i\star_n a),$$
for all $i,j\in \{1,\ldots, n-1\}$ and all $a\in H^*(F_n)$.
Assumption (vi) in Theorem \ref{firstmain} implies that 
$$\sum_{1\le i<j \le n}x_i\star_n x_{j} +\sum_{i=1}^{n-1}q_i=0.
$$
From Theorem \ref{mar} we deduce that $\star_n=\circ$, via the
identifications
$$y_j=\Y_j, \quad q_j =Q_j, \quad 1\le j \le n-1.$$
This implies Theorem \ref{secondmain} for $k=n$:
we also use the fact that the map 
$$\X_i \mapsto \X_{n-i+1}, 1\le i \le n, \quad
Q_j \mapsto Q_{n-j}, 1\le j \le n-1$$
is an automorphism of $qH^*(F_n)$. 
Another consequence of the fact that $\star_n=\circ$ 
concerns the coefficient $\omega_i^{vw}$ in the
expansion given by equation (\ref{yi}):
namely,  the expression
$\omega_i^{vw}|_{q_n=0}$ is prescribed.  

Let us now consider the quotient 
$$(H^*(F_n)\otimes K)/(q_1)
= H^*(F_n)\otimes \bZ[q_2,\ldots, q_{n}]=
\bZ[y_1,\ldots,y_{n-1},q_2,\ldots,q_{n}]/
S(x_1,\ldots,x_n)$$
with the product induced by $\star$. We denote it by $\star_1$.
Assumption  
(vii) in Theorem \ref{firstmain} implies that for any $a\in H^*(F_n)$ and any $d=(0,d_2,\ldots,d_n)$ we have 
\begin{equation}\label{dn1}d_n[(y_1-y_j)\star_1 a]_d=d_j(y_1\star_1 a)_d,
\end{equation}  
as well as
\begin{equation}\label{dn2}d_i[(y_1-y_j)\star_1 a]_d=d_j[(y_1-y_i)\star_1 a]_d,
\end{equation} 
for all $2\le i,j\le n-1$. Let us replace
$$y_1=Y_1,y_1-y_{n-1}=Y_2,y_1-y_{n-2}=Y_3,\ldots,
y_1-y_2=Y_{n-1},$$
and 
$$q_n=Q_1,q_{n-1}=Q_2,\ldots, q_2=Q_{n-1}.$$
As usual, we set $X_j:=Y_j-Y_{j-1}$, for $1 \le j \le n$, where $Y_0=Y_n:=0$.
Then we have
$$X_i=x_{2-i}$$
for all $1\le i \le n$, where the indices are evaluated modulo $n$. 
This implies that
$$H^*(F_n) = \bZ[Y_1,\ldots, Y_{n-1}]/S(X_1,\ldots, X_n).$$
This also implies that 
$$\sum_{1\le i<j \le n}X_i\star_1 X_{j} +\sum_{i=1}^{n-1}Q_i=0.
$$
From equations (\ref{dn1}) and (\ref{dn2}) we deduce that
$$Q_i\frac{\partial}{\partial Q_i}(\Y_j\star_1 a)=Q_j\frac{\partial}{\partial Q_j}
(\Y_i\star_1 a)$$
for all $1\le i,j\le n-1$.
All the hypotheses of Theorem \ref{mar} are verified by $\star_1$.
Consequently, we have  $\star_1=\circ$. 
Thus, for any $1  \le i\le n-1$ and any $v,w\in S_n$, the expression $\omega_i^{vw}|_{q_1=0}$ 
is prescribed. 

A similar reasoning can be made  when we set $q_k$ to zero for an arbitrary $k$.
It is obvious that if $i\in \{1,\ldots,n-1\}$ and  $v,w\in S_n$, there exists {\it at most} one $\omega_i^{vw}\in K$ with all
$\omega_i^{vw}|_{q_1=0},\ldots, \omega_i^{vw}|_{q_n=0}$ prescribed. This proves the uniqueness part of Theorem \ref{firstmain}. Theorem \ref{secondmain} follows also.

\section{Construction of the product $\star$}\label{three}

The construction is related to the periodic Toda lattice of type $A$. The following determinant expansion arises when describing  the conserved
quantities of this integrable system (see the introduction of [Gu-Ot] and
the references therein).
\begin{align}\label{dete}{}& \det \left[ \left(%
\begin{array}{ccccccc}
   Y_1 & q_1 & 0 & \ldots & \ldots & \ldots & -z \\
  -1 & Y_2-Y_1 & q_2 & 0 & \ldots & \ldots & 0 \\
  0 & -1  & Y_3-Y_2 & q_3 & 0 & \ldots & 0 \\
  \hdots & \vdots & \vdots & \vdots & \vdots & \vdots & \hdots \\
 0 & \ldots & 0 & -1 & Y_{n-2}-Y_{n-3} & q_{n-2} & 0 \\
0 & \ldots & \ldots & 0 & -1 & Y_{n-1}-Y_{n-2} & q_{n-1} \\
 q_n/z & \ldots & \ldots & \ldots & 0 & -1 & -Y_{n-1}\\
\end{array}%
\right) +\mu I_n \right]   \\ {}& = \sum_{i=0}^n \R_{i-1}\mu^{n-i} + \R_n \frac{1}{z} - z.\nonumber
\end{align}
We can easily see that  $\R_{-1} = 1$, $\R_0=0$, $$\R_1=\sum_{1\le i<j\le n}(Y_i-Y_{i-1})(Y_{j}-Y_{j-1}) +\sum_{j=1}^nq_j$$ (where
$Y_0=Y_n:=0$),   and $$\R_n=(-1)^{n-1}q_1\ldots q_n.$$
We consider the variables $t_1,\ldots,t_{n-1},q_n$ and  the differential operators
$$ \p^n_i := \frac{\p}{\p t_i} - q_n\frac{\p}{\p q_n},$$ $1\le i\le n-1$.
Then we consider
the ring $$\D^h:=\bC[t_1,\ldots,t_{n-1},e^{t_1}, \ldots, e^{t_{n-1}},q_n, h\p^n_1, \ldots ,h \p^n_{n-1},  h]=\bC[\{t_i\},\{e^{t_i}\},q_n,\{h\p_i^n\},h],$$
where $h$ is a formal variable. This ring is obviously non-commutative.  
It is graded  with respect to 
$$\deg t_i=0, \quad \deg \p^n_i=0, \quad \deg h=2,
\quad   \deg e^{t_i}=4,\quad \deg q_n=4,$$
for $1\le i \le n-1$.  The elements of $\D^h$ should be regarded as differential operators on the space $\bC[t_1,\ldots, t_{n-1},
e^{t_1},\ldots, e^{t_{n-1}}, q_n]$.
In each $\R_k$ we formally replace as follows:
\begin{itemize}
\item $Y_i$ by
$h\p^n_i$, $1\le i\le n-1$,
\item $q_i$ by $e^{t_i}$, for $1\le i \le n-1$ ($q_n$ remains unchanged).
\end{itemize}
The resulting element of $\D^h$ will be denoted by $D_k^h$. 
By considering the expansion of the determinant given by equation
(\ref{dete}) we can see that each ${D}^h_k$ is a linear
combination of  monomials of the form
$(h {\p}/{\p t_1})^{a_1} \ldots (h {\p}/{\p t_{n-1}})^{a_{n-1}}
(hq_n {\p}/{\p q_n})^{a_n}(e^{t_1})^{b_1}\ldots (e^{t_{n-1}})^{b_{n-1}}
q_n^{b_n}$ where
$a_kb_k=0$ for all $1\le k \le n$. This means that 
$D_k^h$ can be obtained by simply making the replacements 
indicated above in the determinant given by (\ref{dete})
and then expanding the determinant 
without any concern about the lack of commutativity of the ring $\D^h$. 

We also consider the quotient ring \begin{equation}\label{first} M^h:=\D^h/\langle D_1^h,\ldots, D_{n-1}^h,
e^{t_1}\cdots e^{t_{n-1}}q_n\rangle,\end{equation}
where the brackets $\langle \ \rangle$ denote the left ideal. We mention that
$$\p_i^n(e^{t_1}\cdots e^{t_{n-1}}q_n)=0,\ \forall i\in\{1,\ldots,n-1\},$$
thus $e^{t_1}\cdots e^{t_{n-1}}q_n$ is a central element of $\D^h$.

We will need the following theorem, which is a straightforward consequence of 
\cite[Lemma 3.5]{Go-Wa}.
\begin{theorem}\label{goodw} {\rm (\cite {Go-Wa})} 
 The ring $$\bC[e^{t_1},\ldots,e^{t_{n-1}},q_n, h\p_1^n,\ldots,h\p_{n-1}^n,h]/
\langle D_1^h,\ldots, D_{n-1}^h\rangle$$
regarded as a left  $\bC[\{e^{t_i}\},q_n,h]$-module has dimension $n!$.
More precisely, any $\bZ$-linear basis $\{ [{f}_w] \ | \ w\in S_n\}$ of 
$H^*(F_n)$, with $f_w\in \bZ[Y_1,\ldots,Y_{n-1}]$,  induces the basis of the previously
mentioned module  which consists of  the cosets of
${f}_w( h\p_1,\ldots, h\p_{n-1} )$, $w\in S_n$.
\end{theorem}

Let us fix a $\bZ$-linear basis $\{[f_w] \ |  \ w\in S_n\}$ of $H^*(F_n)$ as in the previous theorem.
We may assume that $f_{\tau_{i,i+1}}=Y_i$, for all $1\le i\le n-1$, where $\tau_{i,i+1}$ is the transposition 
of $i$ and $i+1$ in $S_n$ (for instance, we can take the 
Schubert basis).
We choose an ordering of $S_n$ (i.e. of the elements of the basis) such that
if  $l(v)< l(w)$ then $v$ comes before $w$.  In this way, the basis
$\{[f_w] \ | \ w\in S_n\}$  of $H^*(F_n)$ consists of $s_0=1$
elements of degree $0$, followed by $s_1$ elements of degree $2$,
$\ldots$, followed by $s_m=1$ elements of degree $2m:=\dim F_n$.
Linear
endomorphisms of $H^*(F_n)\otimes \bC$ can be identified with matrices of size $n!\times n!$ whose entries are in $\bC$.  Any
such matrix  appears as a block matrix of the type
$A=(A_{\alpha\beta})_{1\le\alpha,\beta \le m}$. We  say that a
block matrix $A=(A_{\alpha\beta})_{1\le\alpha,\beta \le m}$ is
$r$-{\it triangular} if $A_{\alpha\beta}=0$ for all $\alpha,\beta$
with $\beta -\alpha <r$. We will use the  splitting of a block matrix
$A=(A_{\alpha\beta})_{1\le \alpha, \beta\le m}$ as
\begin{equation}\label{block}A=A^{[-m]} +\ldots +A^{[-1]}+A^{[0]}+A^{[1]} +\ldots +A^{[m]}
\end{equation}
where each block matrix $A^{[r]}$ is $r$-diagonal.

We set $$P_w:=f_w(h\p_1^n,\ldots, h\p_{n-1}^n),$$   where $w\in S_n$.
By Theorem \ref{goodw}, there exist uniquely determined polynomials
$(\Omega^h_i)_{vw}\in \bC[e^{t_1},\ldots ,e^{t_{n-1}},q_n,h]$ such that
\begin{equation}\label{hpi}h\p_i P_w= \sum_{v\in S_n} (\Omega_i^h)_{vw} P_v \ {\rm modulo} \ 
 \langle D_1^h,\ldots, D_{n-1}^h \rangle,\end{equation}
 for all $1\le i \le n-1$ and all $w\in S_n$. 
 We note that each $(\Omega_i^h)_{vw}$ is a homogeneous polynomial  with respect to the grading given by 
$$\deg h=2, \quad \deg e^{t_i} = 4,1\le i \le n-1, \quad \deg q_n=4 .$$
We consider the matrices
$$\Omega_i^h:=((\Omega_i^h)_{vw})_{v,w\in S_n},\quad \omega_i:=\Omega_i^h \ {\rm
modulo} \ h. $$
Let us also denote by
$\bC[\{t_i\}, \{e^{t_i}\}, q_n]^{(1)}$  the subspace of
$\bC[\{t_i\}, \{e^{t_i}\}, q_n]$ spanned by  the non-constant monomials of  type  $t_1^{k_1}\cdots t_{n-1}^{k_{n-1}}e^{m_1t_1}\cdots e^{m_{n-1}t_{n-1}}q_n^{m_n}$
with the property that if $k_1+\ldots+k_{n-1}\ge 1$ then $m_1\cdots m_n \neq 0$.
The following result plays a  central role in this section.
\begin{proposition}\label{matrix} There exists a  matrix $U=(U_{vw})_{v,w\in S_n}$ of the form
\begin{equation}\label{}U=(I+hV_1+h^2 V_2 +\ldots +h^{m-2}V_{m-2})V_0,\end{equation} with the following properties
\begin{itemize}
\item[(a)]  $V_0$ is a block matrix whose diagonal is $I$,
such that $V_0-I$ is $2$-triangular  and all entries of $V_0$
which are not on the diagonal are   elements of 
$\bC[t_1,\ldots,t_{n-1}, e^{t_1},\ldots,e^{t_{n-1}}, q_n]^{(1)}$
\item[(b)] for any $1\le j\le m-2$, the block matrix $V_j$ is $(j+2)$-triangular and
its entries are elements of 
$\bC[t_1,\ldots,t_{n-1}, e^{t_1},\ldots,e^{t_{n-1}}, q_n]^{(1)}$
\item[(c)] each entry $U_{vw}$ is a  polynomial in
$t_1,\ldots,t_{n-1}, e^{t_1}, \ldots,e^{t_{n-1}}, q_n, h$, homogeneous with respect to $\deg t_i=0, \deg h = 2, \deg e^{t_i}=\deg q_n=4$, such that
$$\deg (U_{vw})+\deg P_v = \deg P_w,$$
\item[(d)] we have \begin{equation}\label{uv} \Omega_i^h U+h\p_i^{n}(U) =UV_0^{-1}\omega_i V_0,\end{equation}
for all $1\le i\le n-1$. 
\end{itemize}
\end{proposition}

In the following we will show how to construct the product $\star$ by using this proposition (its proof will be done after that). For any $w\in S_n$, set
$$\hat{P}_w :=\sum_{v\in S_n} U_{vw} P_v,$$  which is an element of  $\D^h$.
Its coset modulo $\langle D_1^h,\ldots, D_{n-1}^h, e^{t_1}\cdots e^{t_{n-1}}q_n \rangle$ is   \begin{equation}\label{bara}[\hat{P}_w] :=\sum_{v\in S_n} \bar{U}_{vw} [P_v],
 \end{equation}
 where $$\bar{U}_{vw}:=U_{vw} \ {\rm mod}  \ e^{t_1}\cdots e^{t_{n-1}}q_n.$$
Due to  Proposition \ref{matrix}, $\bar{U}_{vw}$ is actually  in $\bC[e^{t_1},\ldots ,e^{t_{n-1}},q_n,h]$
 (i.e. it is free of $t_1,\ldots,t_{n-1}$). 
 On the other hand, equation (\ref{hpi}) implies that in $M^h$ we have
$$ h\p_i^n [P_w]= \sum_{v\in S_n} (\bar{\Omega}_i^h)_{vw} [P_v],$$
where we have denoted 
$$(\bar{\Omega}_i^h)_{vw}:=(\Omega_i^h)_{vw} \ {\rm mod} \ 
e^{t_1}\cdots e^{t_{n-1}}q_n.$$ 
  From equation (\ref{bara}) we deduce
 \begin{equation}\label{flat}h\p_i^n [\hat{P}_w]= \sum_{v\in S_n} (\hat{\Omega}_i^h)_{vw}[\hat{P}_v],
 \end{equation}
 where the matrix $\hat{\Omega}_i^h:=((\hat{\Omega}_i^h)_{vw})_{v,w\in S_n}$ 
 is given by
 \begin{equation}\label{gauge}\hat{\Omega}_i^h=\bar{U}^{-1}\bar{\Omega}_i^h \bar{U} +h\bar{U}^{-1}\p_i^n(\bar{U}).\end{equation}
 Here we have used the elementary fact that if $f\in \bC[t_1,\ldots,t_{n-1},
 e^{t_1},\ldots,e^{t_{n-1}},q_n]$ then we have the equality
 $$\p_i^n\cdot f =\p_i^n(f)+f\cdot \p_i^n$$
 as differential operators. 
Now from equation (\ref{uv}) we deduce that
\begin{equation}\label{haa}\hat{\Omega}_i^h = \bar{V}_0^{-1}\bar{\omega}_i \bar{V}_0,
\end{equation}
 where $$\bar{V}_0:=V_0 \ {\rm modulo} \ e^{t_1}\cdots e^{t_{n-1}}q_n,$$
 and $$\bar{\omega}_i:=\omega_i \ {\rm modulo} \ e^{t_1}\cdots e^{t_{n-1}}q_n.$$
 Indeed, on the one hand, (\ref{uv}) implies that 
 $$(\Omega_i^h \ {\rm mod} )(U \ {\rm mod } ) 
 +h(\p_i^n( U))  {\rm  \ mod }  
 =(U {\rm  \ mod } ) (V_0 {\rm  \ mod } )^{-1}(\omega
 {\rm  \ mod } )( V_0 {\rm  \ mod  } ),$$
 where we have used the abbreviation $${\rm mod} :={\rm modulo } \ e^{t_1}\cdots e^{t_{n-1}}q_n;$$
 on the other hand, we have
 $$\p_i^n (U)  {\rm  \ mod}  
 =\p_i^n (U  {\rm  \ mod } )
 =\p_i^n( \bar{U}).$$  

The cosets $[f_w(Y_1,\ldots,Y_{n-1})]_q$, $w\in S_n$, are a basis of  
\begin{equation}\label{br}\bC[Y_1,\ldots,Y_{n-1},q_1,\ldots,q_n]/( \R_1,\ldots,\R_{n})\end{equation}
 over $K\otimes \bC$. Here we have denoted by $[ \  ]_q$  the cosets modulo
$ ( \R_1,\ldots,\R_{n})$.  Another basis is given by the cosets of 
 \begin{equation}\label{hatf}\hat{f}_w:=\sum_{v\in S_n}(\bar{V}_0)_{vw}f_v,\end{equation}
 $w\in S_n$ (this follows from Proposition \ref{matrix} (a)). In
 $(\bar{V}_0)_{vw}$ we have replaced 
 $e^{t_1},e^{t_2}, \ldots, e^{t_{n-1}}$ by $q_1, q_2, \ldots$, respectively $q_{n-1}$.
 The map
 $$\Phi:H^*(F_n)\otimes K \otimes \bC \to \bC[Y_1,\ldots,Y_{n-1},q_1,\ldots,q_n]/( \R_1,\ldots,\R_{n})$$
 given by $[f_w]\mapsto [\hat{f}_w]_q$, for all $w\in S_n$, is an isomorphism of $K\otimes \bC$-modules.
 We define a product $\star$ on $H^*(F_n)\otimes K\otimes \bC$  by
 $$[f_v]\star[f_w]:=\Phi^{-1}([\hat{f}_v\hat{f}_w]_q),$$
 for all $v,w\in S_n$.
  \begin{proposition}\label{rings}  
 (a) The ring $(H^*(F_n)\otimes K\otimes \bC,\star)$ constructed above satisfies  conditions
 (i)-(vii) in Theorem \ref{firstmain} with $K$ replaced by $K\otimes \bC$.
 
(b) The assertion stated in Theorem \ref{secondmain} with $K$ replaced by $K\otimes 
 \bC$ holds true.

 (c) The subspace $H^*(F_n)\otimes K$ of $H^*(F_n)\otimes K\otimes \bC$ is closed under $\star$. The ring $(H^*(F_n)\otimes K,\star)$ satisfies  conditions
 (i)-(vii) of Theorem \ref{firstmain}. 
 
 (d) 
The map $\Phi$ induces by restriction an isomorphism
of $K$-modules $$\Phi: (H^*(F_n)\otimes K,\star)\to \bZ[Y_1,\ldots,Y_{n-1},q_1,\ldots,q_n]
 /( \R_1,\ldots,\R_{n}).$$
 This is the homomorphism
 of $K$-algebras which maps
  $y_i$ to the coset of $Y_i$, for all $1\le i \le n-1$.  
 \end{proposition}
 \begin{proof} 
 (a)  Only condition (vii) remains to be checked. First we note  that the matrix of multiplication by $Y_i$ on the space described by equation (\ref{br}) with respect to the 
 $K$-basis
  $\{[\hat{f}_w(Y_1,\ldots,Y_{n-1})]_q\}_{w\in S_n}$ is $\hat{\Omega}^h_i|_{e^{t_1}=q_1,\ldots,
  e^{t_{n-1}}=q_{n-1}} \ {\rm modulo} \ q_1 \cdots q_n$ (this follows from
  equations (\ref{haa}), (\ref{hatf})  and the fact that the matrix of multiplication by $Y_i$ with respect to the $K$-basis $\{[f_w(Y_1,\ldots,Y_{n-1})]_q \ | \  w\in S_n\}$ is $\bar{\omega}_i$). 
  From equation (\ref{flat}), the identity $\p_i^n\p_j^n=\p_j^n\p_i^n$ for all $1\le i,j\le n-1$, and the fact that all $\hat{\Omega}_i^h$ are independent of $h$, we deduce that each entry of the matrix $\p_i^n(\hat{\Omega}_j^h)-\p_j^n(\hat{\Omega}_i^h)$ is a multiple of
  $e^{t_1}\cdots e^{t_{n-1}}q_n$, thus it is equal to zero in $K$.  
 This implies condition (vii). 

(b) The assertion can be proved by using the same method as in Section 2.

(c) The assertion follows immediately from the following two claims.

\noindent{\it Claim 1.} For any $i\in \{1,\ldots,n-1\}$ and any $a\in H^*(F_n)$, the product
$y_i\star a$ is in $H^*(F_n)\otimes K$.

\noindent{\it Claim 2.} Any $a\in H^*(F_n)$ can be written as
$a=f(\{y_i\star\},\{q_j\})$, where $f$ is a polynomial in $\bZ[\{Y_i\},\{q_j\}]$.

To prove Claim 1, we consider again the Schubert basis
$\{\sigma_w \ | \ w\in S_n\}$ of $H^*(F_n)$.
We take $i\in\{1,\ldots,n-1\}$ and $w\in S_n$ and  consider
$\omega_i^{vw}\in K\otimes \bC$ given by 
$$y_i\star \sigma_w = \sum_{v\in S_n}\omega_i^{vw}\sigma_v.$$
We show that $\omega_i^{vw}\in K$. 
This is an immediate consequence of the fact that
for any $k\in\{1,\ldots,n\}$ we have
$$\omega_i^{vw}|_{q_k=0}\in \bZ[q_1,\ldots,q_{k-1},q_{k+1},\ldots,q_n].$$
To prove this,  we take into account that, by point (b) above, the ring 
$((H^*(F_n)\otimes K \otimes \bC)/(q_k),\star)$ is isomorphic to the quantum cohomology ring of $F_n$ tensored with $\bC$. Moreover,                  
$\omega_i^{vw}|_{q_k=0}$
is a coefficient of the expansion of the quantum product of two elements of $H^*(F_n)$
(one of them of degree 2).

To prove Claim 2, we assume that $a$ is homogeneous and we use induction by
$\deg a$. If $a$ has degree 0 or 2, the claim is obvious.
At the induction step, we assume that the claim is true for any element of $H^*(F_n)$
of degree at most $d$  and show that it is true for  any $y_ia$, where 
$i\in \{1,\ldots,n-1\}$ and 
$a\in H^*(F_n)$ with $\deg a \le d$. To this end, it is sufficient to note
that $y_i \star a- y_i a$ is of the form described in the claim: this follows from Claim 1,  point a)
(especially conditions (ii) and (iii) in Theorem \ref{firstmain}), and the induction hypothesis. 
 
(d) We use the definition of $\star$ on $H^*(F_n)\otimes K \otimes \bC$ and  point c)
above.

\end{proof} 

Now comes the postponed proof of Proposition \ref{matrix}. 

\noindent {\it Proof of Proposition \ref{matrix}.} 
We go along the lines of Amarzaya and Guest [Am-Gu, Section 2]. 
Namely, for each $i\in\{1,\ldots,n-1\}$ we write
$$\Omega_i^h =\omega_i+h\theta^{(1)}_i +h^2\theta^{(2)}_i
+\ldots +h^{m-2}\theta^{(m-2)}_i,$$
where the matrix $\omega_i$ is $-1$-triangular and $\theta^{(k)}_i$ is 
$k+1$-triangular, $1\le k\le m-2$. 
We have the following splittings (see (\ref{block})):
\begin{align*}V_0=&I+V_0^{[2]}+V_0^{[3]}+\ldots + V_0^{[m]}\\
V_{\ell}=&V_{\ell}^{[\ell+2]} +V_{\ell}^{[\ell+3]}+\ldots + V_{\ell}^{[m]} \ \ (1\le \ell \le m-2)\\
\omega_i=&\omega_i^{[-1]} +\omega_i^{[0]} +\omega_i^{[1]}+\ldots
+\omega_i^{[m]}\\
\theta^{(k)}_i=&\theta^{(k),[k+1]}_i+\theta^{(k),[k+2]}_i+\ldots
+\theta^{(k),[m]}_i \ \ (1\le k\le m-1, 1\le i \le n-1).
\end{align*}
The main point is that there exists a   total ordering on the matrices
$V_{\ell}^{[r]}$, $r\ge \ell+2\ge 2$ in such a way that equation (\ref{uv}) is equivalent to the system\begin{equation}\label{sys5}\p_i^n(V_{\ell}^{[r]}) = 
{\rm \ expression \ involving }\ 
V_{\ell'}^{[r']} >
V_{\ell}^{[r]}, \ \forall i\in\{1,\ldots,n-1\}\end{equation}  for  $r\ge \ell+2\ge 2$.
In this way we  can determine all $V_{\ell}^{[r]}$ inductively. 

 The system is compatible, the compatibility condition being
$$\p_i^n(\Omega_j^h)-\p_j^n(\Omega_i^h)=\frac{1}{h}[\Omega_j^h,\Omega_i^h],$$ where $1\le i,j\le n-1$, which in turn, follows easily from (\ref{hpi})
and the fact that $\p_i^n\p_j^n=\p_j^n\p_i^n$. It is also important to note that the right hand side of equation (\ref{sys5}) is a linear combination of $\theta^{(k)}_i$ and left or right side products of such matrices   with one or two  $V_{\ell'}^{[j']}$
(see the equation displayed before Definition 2.3 in [Am-Gu]). By using this we show recursively that
all $V_{\ell}^{[j]}$ have  entries in $\bC[t_1,\ldots,t_{n-1}, e^{t_1},\ldots,e^{t_{n-1}}, q_n]^{(1)}$. This follows from the general fact that if $f\in
\bC[t_1,\ldots,t_{n-1}, e^{t_1},\ldots,e^{t_{n-1}}, q_n]$
has the property that
$$\p_i(f)\in  \bC[t_1,\ldots,t_{n-1}, e^{t_1},\ldots,e^{t_{n-1}}, q_n]^{(1)},\ 
{\rm for \ all \ } 1\le i \le n-1,$$ then $f \in  \bC[t_1,\ldots,t_{n-1}, e^{t_1},\ldots,e^{t_{n-1}}, q_n]^{(1)}$. The proposition follows in an elementary way.
\hfill $\square$

\section{Frobenius property and orthogonality}

We first  prove Theorem \ref{onemain}.

\noindent{\it Proof of Theorem \ref{onemain}.} We consider for each $i\in \{1,\ldots,n-1\}$ the 
 $K$-linear operator $\fY_i$ on $H^*(F_n)\otimes K$ given by
 $$(\fY_i(a),b)=(y_i\star b,a),$$
 for all $a,b\in H^*(F_n)$.

Now let us consider again the Schubert basis
$\{\sigma_w \ | \ w\in S_n\}$ of $H^*(F_n)$. We know that it satisfies the
orhogonality condition
\begin{equation}\label{orthogonal}(\sigma_v,\sigma_w)=\begin{cases}
0,  \ {\rm if} \  v\neq w_0w\\
1, \ {\rm if} \ v=w_0w
\end{cases}
\end{equation}
where $w_0$ is the longest element of $S_n$. 
(That is, $w_0(i)=n-i+1$, for all $i\in \{1,\ldots, n\}$.) 
We deduce easily that the matrix of 
$\fY_i$  on $H^*(F_n)\otimes K$ with respect to the basis $\{\sigma_{w_0w} \ | \ w\in S_n\}$  is the transposed of the matrix of  ``$y_i\star$" with respect to the basis 
$(\sigma_{w})_{w\in S_n}$. From this we deduce as follows:
\begin{itemize}
\item[1.] for any $i,j\in \{1,\ldots,n-1\}$,
the operators $\fY_i$ and $\fY_j$ commute with each other,
\item[2.] if $a$ is a homogeneous element of $H^*(F_n)$, then $\deg\fY_i(a)=\deg a +2$
and  $\fY_i(a)$ modulo $(q_1,\ldots,q_n)$ is equal to $y_ia$.
\end{itemize}

We will use the following claim.

\noindent{\it Claim 1.} For any $a\in H^*(F_n)$ there exists a polynomial
$f\in \bZ[\{Y_i\},\{q_j\}]$ with $$a=f(\{\fY_i\}, \{q_j\}).1.$$

To prove this, we  use the same method as for  Claim 1 in the proof of Proposition 
\ref{rings} (property 2   above is used here).

 This allows us to define a product, denoted by $\bullet$, on $H^*(F_n)\otimes K$, as follows.
 For any $f,g\in \bZ[\{Y_i\}, \{q_j\}]$, we set
 $$\left(f(\{\fY_i\}, \{q_j\}).1\right) \bullet \left(g(\{\fY_i\}, \{q_j\}).1\right)
 :=(fg)(\{\fY_i\}, \{q_j\}).1,$$
 where ``$.1$" denotes evaluation on the element 1 of $H^*(F_n)$.
 This product satisfies the hypotheses (i)-(vii) of Theorem \ref{firstmain}. 
 The first non-obvious one is (vi). This follows from the fact that
 $$y_i\bullet y_j =\fY_i\fY_j .1,$$ which implies that
 $$(y_i\bullet y_j,a)=(\fY_i\fY_j. 1,a)= (y_i\star a, \fY_j. 1)=(y_j\star y_i\star a,1),$$
 for all $a\in H^*(F_n)$. We use that $\star$ is associative and  satisfies condition (vi). 

Let us now check condition (vii). Because $\fY_i.1=y_i$, from  Claim 1 we deduce that  
\begin{equation}\label{bullet}y_i\bullet a =\fY_i(a),\end{equation}
for all $i\in \{1,\ldots, n-1\}$ and all $a\in H^*(F_n)$.   We use again the fact that the matrix of 
$\fY_i$  on $H^*(F_n)\otimes K$ with respect to the basis $\{ \sigma_{w_0w}
\ | \ w\in S_n\}$  is the transposed of the matrix of  ``$y_i\star$" with respect to the basis 
$(\sigma_{w})_{w\in S_n}$. Condition (vii) follows. 
From Theorem \ref{firstmain}, we deduce that $\bullet$ is the same as $\star$.

\noindent {\it Claim 2.} We have 
$$y_{i_1}\star \ldots \star y_{i_k} =\fY_{i_1} \ldots  \fY_{i_k}.1,$$
for all $k\ge 1$ and all $i_1,\ldots,i_k\in \{1,\ldots,n-1\}$. 

We prove the claim by induction on $k$. For $k=1$ the statement is clear.
Assume that it is true for $k$. For any $i_1,\ldots,i_{k+1}\in \{1,\ldots,n-1\}$ we have
$$\fY_{i_1} \ldots  \fY_{i_{k+1}}.1=y_{i_1}\bullet (\fY_{i_2} \ldots  \fY_{i_{k+1}}.1) 
=y_{i_1}\star (y_{i_2} \star \ldots \star y_{i_{k+1}}),$$
where we have used the induction hypothesis and equation (\ref{bullet}). 

Since the ring $(H^*(F_n)\otimes K,\star)$ is generated as a $K$-algebra by $y_1,\ldots,y_{n-1}$, Claim 2  implies immediately that
$$\fY_i(a) =y_i \star a$$
for all $i\in \{1,\ldots, n-1\}$ and all $a\in H^*(F_n)$. 
From the definition of $\fY_i$ we deduce that
$$(y_i\star b,c) =(y_i\star c,b)$$
for all $i\in \{1,\ldots,n-1\}$ and all $b,c\in H^*(F_n)$. 
In turn, this implies
$$(a\star b,c) =(a\star c,b)$$ 
for all $a,b,c\in H^*(F_n)$. 
Theorem \ref{onemain} is now proved.\hfill $\square$

In the second part of this section we consider again the presentation of
$(H^*(F_n)\otimes K,\star)$ in terms of generators and relations  given in the previous section. More precisely, let $\Phi$ be the isomorphism  described by Proposition
\ref{rings} (d).
For each $w\in S_n$ we pick a polynomial $f_w\in  \bZ[\{Y_i\},\{q_j\}]$ such that
$$\Phi(\sigma_w) =[f_w]_q,$$
where $[ \ ]_q$ denotes cosets of polynomials modulo  $( \R_1,\ldots,
 \R_n) $ (like in the previous section). For any $f\in \bZ[\{Y_i\},\{q_j\}]$ we consider the decomposition 
$$[f]_q= \sum_{w\in S_n} \alpha_w[f_w]_q$$
where $\alpha_w \in K$. Then we set
$$\la [f]_q\ra:=\alpha_{w_0}$$
where $w_0$ is the longest element of $S_n$ (see above).
If $f$ and $g$ are in $\bZ[\{Y_i\},\{q_j\}]$, we define
$$\la[f]_q,[g]_q\ra:=\la [fg]_q\ra.$$
We show that the polynomials $f_w$, $w\in S_n$, satisfy the following orthogonality relation, similarly to   the quantum Schubert polynomials (cf. \cite{Fo-Ge-Po} and
\cite{Ki-Ma}).

\begin{proposition}\label{lang} We have
$$\la [f_v]_q,[f_w]_q\ra =\begin{cases}
0,  \ {\rm if} \  v\neq w_0w\\
1, \ {\rm if} \ v=w_0w.
\end{cases}
$$
\end{proposition}

\begin{proof} We consider the expansion
$$[f_vf_w]_q=\sum_{u\in W}\alpha_u[f_u]_q.$$
Since $\Phi$ is a ring isomorphism, this implies that
$$\sigma_v\star\sigma_w = \sum_{u\in W}\alpha_u\sigma_u.$$
We deduce that
$$\alpha_{w_0}=(\sigma_v\star\sigma_w,1)=(\sigma_v,\sigma_w),$$
where we have used the Frobenius property. 
Finally, we use the orthogonality relation (\ref{orthogonal}).
\end{proof}
\bibliographystyle{abbrv}

\begin{thebibliography}{Fo-Ge-Po}

\bibitem[Am-Gu]{Am-Gu} A. Amarzaya and M.  A.  Guest, {\em Gromov-Witten invariants of flag
manifolds, via D-modules}, Jour. London Math. Soc. (2) {\bf 72} (2005), 121--136


\bibitem[Fo-Ge-Po]{Fo-Ge-Po}  S. Fomin, S. Gelfand, and A. Postnikov, {\em Quantum
Schubert polynomials},  J. Amer. Math. Soc.,  {\bf 10} (1997),
565--596

\bibitem[Fu]{Fu} W. Fulton,
{\it Young Tableaux}, London Math. Soc. Student Texts, vol. 35, Cambridge
University Press, Cambridge 1997 

\bibitem[Fu-Pa]{Fu-Pa} W. Fulton and R. Pandharipande,
{\em Notes on stable maps and quantum cohomology},  Algebraic 
Geometry - Santa Cruz 1995, Proc. Sympos. Pure Math., 62, Part 2, editors J. Kollar, R. Lazarsfeld, 
and D.R. Morrison, 1997, 45-96

\bibitem[Gi-Ki]{Gi-Ki} A. Givental and B. Kim,
{\it  Quantum cohomology of flag manifolds and Toda lattices},
Comm. Math. Phys. {\bf 186} (1995), 609--641


\bibitem[Go-Wa]{Go-Wa} R. Goodman, N.R. Wallach, {\it Classical and quantum-mechanical systems of Toda lattice 
type} III, Comm. Math. Phys. {\bf 105} (1986), 473-509 

\bibitem[Gu-Ot]{Gu-Ot} M. Guest and T. Otofuji, {\it Quantum cohomology and the periodic Toda lattice}, Comm. 
Math. Phys. {\bf 217} (2001), 475--487



\bibitem[Ki-Ma]{Ki-Ma} A. Kirillov and T. Maeno,
{\it Quantum double Schubert polynomials, quantum Schubert polynomials and Vafa-Intriligator formula}, Discrete Math. {\bf 217} (2000), 191-224

\bibitem[Ko]{Ko} B. Kostant,
{\it The solution to a generalized Toda lattice and representation theory}, Adv. Math. 
{\bf 34} (1979), 195--338 


\bibitem[Ma1]{Ma1} A.-L. Mare, {\it Quantum cohomology of the infinite dimensional flag manifolds}, Adv. Math. {\bf 185} (2004), 347--369

\bibitem[Ma2]{Ma2} A.-L. Mare, {\it A characterization of the quantum cohomology ring of  $G/B$ and applications}, Canad. Jour. Math. {\bf 60} (2008),  875--891


\end{thebibliography}

\section{Polynomial representatives of Schubert classes: an example}

Let us consider again the presentation of 
$(H^*(F_n)\otimes K,\star)$ given in Proposition \ref{rings} (d).
As usual in this paper, we set
$$X_i:=Y_i-Y_{i-1}, 1\le i\le n, \ {\rm where} \   Y_0=Y_n:=0.$$
In this way $\R_1,\ldots,\R_{n-1}$ defined at the beginning of Section \ref{three}
are polynomials in $X_1,\ldots,X_n$; we set
$$\R_0:=X_1+\ldots +X_n.$$
We have the ring isomorphism
\begin{equation}\label{schub}(H^*(F_n)\otimes K,\star)\simeq
\bZ[X_1,\ldots,X_n,q_1,\ldots,q_n]/( \R_0,\R_1,\ldots,\R_{n}),\end{equation}
 where the classes $x_1:=y_1,x_2:=y_2-y_{1},\ldots,x_{n-1}:=y_{n-1}-y_{n-2},
x_n:=-y_{n-1}$ correspond to the cosets of $X_1,X_2, \ldots,X_{n-1}$, respectively
$X_n$. 

In this section we are concerned with the problem of finding, for an arbitrary $w\in S_n$, of a  polynomial $f_w\in \bZ[X_1,\ldots,X_n,q_1,\ldots,q_n]$ like in the previous section: its coset corresponds to the Schubert  class $\sigma_w$ via
the isomorphism (\ref{schub}). 
In general, this problem can be solved by using Theorem \ref{secondmain} and the 
knowledge of the quantum cohomology ring of $F_n$.
We will illustrate this idea in the special case of $F_3$.

The Schubert classes for $F_3$ are as follows (cf. e.g. \cite[Tables 14]{Fo-Ge-Po}):
\begin{align*}
{}&\S_{123}=1\\
{}& \S_{213}=x_1\\
{}& \S_{132}=x_1+x_2\\
{}& \S_{231}=x_1x_2\\
{}& \S_{312}=x_1^2\\
{}& \S_{321}=x_1^2x_2.
\end{align*}

For each of the  classes $\S_w$  above, we look for a polynomial
$f_w(X_1,X_2,X_3,q_1,q_2,q_3)$ which is homogeneous
relative to $\deg X_i =2,\deg q_j=4$, has the same degree as the cohomology class $\S_w$,  has integer coefficients, and  the element
$f_w(x_1\star,x_2\star,x_3\star,q_1,q_2,q_3)$ of $H^*(F_3)\otimes
\bZ[q_1,q_2,q_3]$ is the same as $\S_w$.
For instance, for $\S_{312}$, we look for
$f=X_1^2+g(q_1,q_2,q_3)$, where $g$ is a homogeneous polynomial of degree 1.
The condition which must be satisfied is
\begin{equation}\label{xcicx}x_1\star x_1 +g(q_1,q_2,q_3)=x_1^2.\end{equation}
We make successively $q_k=0$, for $k=1,2,3$.

First we make $q_3=0$. Equation (\ref{xcicx}) implies an identity in 
$qH^*(F_3)$ obtained by making the replacements prescribed by Theorem \ref{secondmain}.
They are as follows:
$$\X_1:=x_1, \X_2:=x_2, \X_3=x_3,Q_1=q_1,Q_2:=q_2.$$
The identity in $qH^*(F_3)$ is
$$\X_1\circ \X_1+g(Q_1,Q_2,0)=\X_1^2.$$
On the other hand, we know that (see for instance [Am-Gu, Section 3]) 
$$\X_1\circ \X_1=\Y_1\circ \Y_1 = \Y_1^2+Q_1,$$
thus
$g(Q_1,Q_2,0)=-Q_1$, which gives 
\begin{equation}\label{una}g(q_1,q_2,0)=-q_1.\end{equation}

We make $q_1=0$. The replacements are as follows:
$$\X_1:=x_1, \X_2:=x_3,\X_3:=x_2, Q_1:=q_3,Q_2:=q_2.$$
The identity in $qH^*(F_3)$ is
$$\X_1\circ \X_1+g(0,Q_2,Q_1)=\X_1^2.$$
Using the same method as above, we deduce
\begin{equation}\label{doua}g(0,q_2,q_3)=-q_3.\end{equation}

Finally, we make $q_2=0$. The replacements are
$$\X_1:=x_2, \X_2:=x_1,\X_3:=x_3, Q_1=q_1,Q_2=q_3.$$
The identity in $qH^*(F_3)$ is
$$\X_2\circ \X_2 +g(Q_1,0,Q_2)=\X_2^2.$$
This time we obtain
\begin{equation}\label{tre}g(q_1,0,q_3) = -q_1-q_3.\end{equation}

From equations (\ref{una}), (\ref{doua}), and (\ref{tre}) we deduce that
$$g(q_1,q_2,q_3)=-q_1-q_3.$$
Thus, a polynomial representative of $\S_{312}$ is 
$$X_1^2-q_1-q_3.$$

Similarly one obtains the desired polynomials for all Schubert classes.
They are as follows:

\begin{align*}
{}& \S_{123} \ : \ 1\\ 
{}& \S_{213} \ :  \  X_1\\
{}& \S_{132} \ :  \  X_1+X_2\\
{}& \S_{231} \ :  \  X_1X_2+q_1\\
{}& \S_{312} \ :  \  X_1^2-q_1-q_3\\
{}& \S_{321} \ : \  X_1^2X_2+q_1X_1-q_3X_2
\end{align*}

\end{document}